\documentclass{amsart}

\usepackage{amsmath}
\usepackage{amsthm}
\usepackage{amssymb}
\usepackage[all]{xy}

\newtheorem{thm}{Theorem}
\newtheorem{lemma}[thm]{Lemma}
\newtheorem{prop}[thm]{Proposition}

\newtheorem{cor}[thm]{Corollary}

\newtheorem{thma}{Theorem}

\theoremstyle{definition}
\newtheorem{defn}[thm]{Definition}
\newtheorem{eg}[thm]{Example}

\theoremstyle{remark}
\newtheorem{rem}[thm]{Remark}
\newtheorem{notn}[thm]{Notation}

\numberwithin{thm}{section}
\numberwithin{equation}{section}

\newcommand {\Fp} {\ensuremath {\mathbb{F}_p} }
\newcommand {\F} {\ensuremath {\mathbb{F}} }
\newcommand {\freeT} {\ensuremath {\mathbb{T}} }

\newcommand {\Z} {\ensuremath {\mathbb{Z}} }
\newcommand {\Q} {\ensuremath {\mathbb{Q}} }

\newcommand {\CP} {\ensuremath {\mathbb{CP}} }
\newcommand {\AH}[1] {\ensuremath {\mathbf{A} (#1)} }

\newcommand {\lX} {\ensuremath {\Omega X} }
\newcommand {\lY} {\ensuremath {\Omega Y} }
\newcommand {\lZ} {\ensuremath {\Omega Z} }

\newcommand {\freeL} {\ensuremath {\mathbb{L}} }
\newcommand {\isom} {\ensuremath {\cong} }
\newcommand {\tensor} {\ensuremath {\otimes} }
\newcommand {\incl} {\ensuremath {\hookrightarrow} }

\newcommand {\sdp} {\ensuremath {\ltimes} }
\newcommand {\PiS} {\ensuremath {\prod \mathcal{S}}}

\newcommand {\bB} {\ensuremath {\mathbf{B}} }
\newcommand {\nP} {\ensuremath {\tilde{P}} }
\newcommand {\fpnP} {\ensuremath {\forall p \in \nP } }

\newcommand {\feachpnP} {for each $p \in \nP$}
\newcommand {\onto} {\ensuremath {\twoheadrightarrow} }
\newcommand {\eL} {\ensuremath {\mathbf{\underline{L}}} }

\newcommand {\isomto} {\ensuremath {\xrightarrow{\isom}} }
\newcommand {\xto}[1] {\ensuremath {\xrightarrow{#1}} }
\newcommand {\eolBox} {\vspace{2.5mm} \ensuremath {\hfill \Box} }

\DeclareMathOperator{\im}{im}
\DeclareMathOperator{\Rank}{Rank}

\DeclareMathOperator{\gr}{gr}

\DeclareMathOperator{\coker}{coker}

\DeclareMathOperator{\Torsion}{Torsion}

\begin{document}

\title{{Free and semi-inert cell attachments}}
\author{Peter Bubenik}
\date \today
\address{D\'epartment de Math\'ematics, Ecole Polytechnique F\'ed\'erale de
Lausanne \\ CH-1015, Lausanne, Switzerland}
\email{peter.bubenik@epfl.ch}
\keywords{cell attachments, loop space, loop space homology,
Adams-Hilton models}
\subjclass[2000]{Primary 55P35 ; Secondary 16E45}

\begin{abstract}
Let $Y$ be the space obtained by attaching a finite-type wedge of cells to 
a simply-connected, finite-type CW-complex.

We introduce the \emph{free} and \emph{semi-inert} conditions on the
attaching map which broadly generalize the previously studied
\emph{inert} condition.  
Under these conditions we determine $H_*(\Omega Y;R)$ as an $R$-module and
as an $R$-algebra respectively. 
Under a further condition we show that $H_*(\Omega Y;R)$ is
generated by Hurewicz images.

As an example we study an infinite family of spaces constructed using
only semi-inert cell attachments.
\end{abstract}

\maketitle

\section {Introduction} \label{sectionIntro}

In this article we will work in the usual category of pointed,
simply-connected topological spaces with the homotopy-type of
finite-type CW-complexes. 
We will assume that the ground ring $R$ is either $\Fp$ with $p>3$ or
is a subring of $\Q$ which contains $\frac{1}{6}$.

We are interested in the following problem, perhaps first studied by
J.H.C. Whitehead around 1940~\cite{whiteheadAddingRelations,
whiteheadSimplicialSpaces}. 

\noindent
\textbf{The cell attachment problem:}
\emph{
Given a topological space $X$, what is the effect on the loop space
homology and the homotopy-type if one attaches one or more cells to
$X$?
}

We approach this problem from the point of view that one is interested
in understanding finite cell complexes localized away from finitely
many primes~\cite{anickSLSD}.

Given a space $X$ and a map $f: W \to X$ where $W = \vee_{j \in J}
S^{n_j}$, the \emph{adjunction space}
\[
Y = X \cup_f \bigl( \bigvee_{j \in J} e^{n_j+1} \bigr),
\]
is a homotopy cofibre of $f$.
Let $i$ denote the inclusion $X \incl Y$.

The cell attachment problem has been studied in two special cases.
One approach is to place a strong condition on the space $X$.
This was done by Anick~\cite{anickCat2} who considered the case
where $X$ is a wedge of spheres.
Another approach is to place a strong condition on the attaching map
$f$.
This was done by Lemaire and Halperin
\cite{lemaireAutopsie},\cite{halperinLemaireInert} and 
F\'{e}lix and Thomas~\cite{felixThomasAttach} who assumed that $f$ is
\emph{inert}.

The attaching map $f:W \to X$ is said to be \emph{inert} over a ring
$R$ if the induced map $H_*(\Omega i; R): H_*(\lX;R) \to H_*(\lY;R)$ is
a surjection.

In this article we generalize these two approaches, with our
development following~\cite{anickCat2}.
We generalize Anick's assumption to the more general condition that
$H_*(\lX;R)$ is $R$-free and is generated by Hurewicz images. 
This is trivial in the case where $R = \Q$, and we will give
conditions under which this holds for more general coefficient rings
(see Corollary~\ref{corFiniteCellComplex}).
Furthermore we give two generalizations of the inert condition, one of
which is strictly stronger than the other.

We now define these conditions in the case where $R$ is a field.
We will use the following notation.

\begin{notn}
Given a space $X$, let $L_X$ denote the image of the Hurewicz map
$h_X: \pi_*(\lX) \tensor R \to H_*(\lX;R)$.
$L_X$ is a graded Lie algebra under the commutator bracket of the
Pontrjagin product.
Given a map $W \to X$, let $L^W_X$ denote the image of the induced Lie
algebra map $L_W \to L_X$. 
Note that the map is omitted from the notation.
Let $[L^W_X] \subset L_X$ denote the Lie ideal generated by $L^W_X$.
\end{notn}

\begin{defn}
Define a cell attachment $f:W \to X$ to be \emph{free} if $[L^W_X]$ is
a free Lie algebra. 
\end{defn}

There are examples of spaces given by non-free cell attachments that
can be constructed by free cell attachments if one changes the order in
which the cells are attached
(eg. $\CP^2 \cup_f e^3$ where $f$ is the inclusion $S^2 \incl \CP^2$,
see~\cite[Example 4.5]{hessLemaireNice}).
When $R$ is a field, it is conceivable that any cell complex can be
constructed using only free cell attachments if one chooses an
appropriate cellular structure. 

Free cell attachments are convenient to work with because of the
following fact about universal enveloping algebras, which we prove as
Lemma~\ref{lemmaI}.
Since it may be of independent interest we state it here.
$U$ denotes the universal enveloping algebra functor.

\begin{prop}
Let $L$ be a connected, finite-type Lie algebra over a field. 
Let $J$ be a Lie ideal of $L$ which is a free Lie algebra, $\freeL
W$.
Take $I$ to be the two-sided ideal of $UL$ generated by $J$.
Then the multiplication maps $UL \tensor W \to I$ and $W \tensor UL \to
I$ are isomorphisms of left and right $UL$-modules respectively.
\end{prop}

Assume that $Y$ is obtained by a free cell attachment $\bigvee_{j
\in J} S^{n_j} \xto{\bigvee \alpha_j} X$.
Let $\hat{\alpha_j}$ denote the adjoint of $\alpha_j$.
We will show that $H_*(\lY;R)$ can be determined by calculating the
homology of the following simple differential graded Lie algebra
(dgL)
\[ \eL = (L_X \amalg \freeL \langle y_j \rangle_{j \in J}, d), 
\text{ where }
dy_j = h_X(\hat{\alpha_j}).
\]
In addition to the usual grading, $\eL$ has a second grading given by
letting $L_X$ be in degree $0$ and letting each $y_j$ be in degree
$1$.
Remarkably, we will show that for free cell attachments one only needs
to calculate $H\eL$ in degrees $0$ and $1$.

Let $(H\eL)_i$ denote the component of $H\eL$ in degree $i$.
For degree reasons, $(H\eL)_0$ acts on $(H\eL)_1$ by the adjoint action.
\begin{defn}
Define $f$ to be a \emph{semi-inert} cell attachment if it is a free
cell attachment and $(H\eL)_1$ is a free $(H\eL)_0$-module. 
\end{defn}

We will see in Section~\ref{sectionAH} that there is an obvious
filtration on $H_*(\lY;R)$. 
Let $\gr_*(H_*(\lY;R))$ be the \emph{associated graded object}.
We will show (see Definition~\ref{defnSemiInertCA})
that a free cell attachment is semi-inert if and only if
$\gr_1(H_*(\lY;R))$ is a free $\gr_0(H_*(\lY;R))$-bimodule. 


\begin{thma} \label{thmA}
Let $Y = X \cup_f \left(\bigvee_{j \in J} e^{n_j+1}\right)$.
Assume that 
$f$ is \emph{free}. \\
(i) Then as algebras 
\[ \gr(H_*(\lY; R)) \isom U( (H\eL)_0 \sdp \freeL  ((H\eL)_1) )
\]
with $(H\eL)_0 \isom {L}_X/[{L}^W_X]$ as Lie algebras. \\ 
(ii) Furthermore if $f$ is \emph{semi-inert} then for some $K'$,
\[ H_*(\lY;R) \isom U ( L^X_Y \amalg \freeL K' )
\]
as algebras.
\end{thma}

The above result is given more precisely in Theorem~\ref{thmPreciseA}.
This theorem is a nearly direct translation of a purely algebraic
result given in Theorem~\ref{thmDGAextension}.
This algebraic result may have other applications such as the
calculation of the mod $p$ Bockstein spectral sequence  (BSS) of finite
CW-complexes.
Scott~\cite{scottBSS} has shown that for sufficiently large $p$ each
term in the mod $p$ BSS of such spaces is a universal enveloping
algebra of a dgL.

\begin{cor} \label{corNice}
Free cell attachments are \emph{nice} in the sense of Hess
and Lemaire \cite{hessLemaireNice}.
\end{cor}

When $R=\Q$, Milnor and Moore~\cite{milnorMooreHopfAlgebras} proved
the fabulous result that the canonical algebra map $U(\pi_*(\lY)
\tensor \Q) \to H_*(\lY;\Q)$ is an isomorphism.
Scott~\cite{scottTFMM} generalized this result to $R \subset \Q$ for
finite CW-complexes when certain primes are invertible in $R$.
For an $R$-module $M$, let $FM = M / \Torsion(M)$.
Let $P$ denote the primitive elements of $FH_*(\lY;R)$ which are a Lie
subalgebra.
Scott showed that $UP \isomto FH_*(\lY;R)$ and that $F(\pi_*(\lY)
\tensor R)$ injects into $P$. 
However he showed that in general this injection is not a surjection.

We will give sufficient conditions 
under which one obtains the desired isomorphism 
$UF(\pi_*(\lY) \tensor R) \to FH_*(\lY;R)$

In Section~\ref{sectionHurewicz} we will assume that the Hurewicz map
$h_X$ has a right inverse. 
Using this map we will define the set of \emph{implicit primes} of $Y$. 
Intuitively, they are the primes $p$ for which $p$-torsion is used in the
attaching map $f$.
Let 
\[ \mathcal{S} = \{ S^{2m-1},\,\Omega S^{2m+1} \:|\: m \geq 1\}.
\]
Let $\PiS$ be the collection of spaces homotopy equivalent
to a weak product of spaces in $\mathcal{S}$.

\begin{thma} \label{thmB}
Let $Y = X \cup_f \left(\bigvee e^{n_j+1}\right)$.
Assume that $f$ is free, that the Hurewicz map $h_X$ has a right
inverse, and that the implicit primes are invertible.
Then the canonical algebra map
\begin{equation} \label{eqnUL_Y} 
UL_Y \to H_*(\lY;R)
\end{equation} 
is a surjection.
Furthermore if $R \subset \Q$ then \eqref{eqnUL_Y} is an isomorphism
and localized at $R$, $\lY \in \PiS$.
If in addition $f$ is semi-inert, then 
\[ L_Y \isom H\eL \isom (H\eL)_0 \amalg \freeL ((H\eL)_1)
\]
as Lie algebras, and $h_Y$ has a right inverse.
\end{thma}

Note that the surjection of~\eqref{eqnUL_Y} implies that $H_*(\lY;R)$
is generated as an algebra by Hurewicz images. 
Again, more details are given in Theorem~\ref{thmHurewicz}.

\begin{cor} \label{corTFMM}
If $R \subset \Q$ then the canonical algebra map 
\[ UF(\pi_*(\lY) \tensor R) \to H_*(\lY;R)
\]
is an isomorphism.
\end{cor}

\begin{cor} \label{corFiniteCellComplex}
If $Z$ is a finite cell complex constructed using only
semi-inert cell attachments then localized away from a finite set of primes,
$\lZ \in \PiS$.
\end{cor}

It is a long-standing conjecture of
Avramov~\cite{avramovFreeLieSubalgebras} and F\'elix~\cite{fhtCat2}
that if $Z$ has {finite LS category} then $L_Z$ is either finite
dimensional or contains a free Lie subalgebra on two generators.
Our final corollary provides further support for this conjecture.

\begin{cor} \label{corFelixAvramov}
If $R \subset \Q$, $f$ is a free cell attachment and $\dim(H\eL)_1 >
1$ then $L_Y$ contains a free Lie subalgebra on two generators. 
\end{cor}

We conclude by giving examples of spaces constructed
out of semi-inert cell attachments, together with their Hurewicz
images.
In particular, we give an infinite family of finite CW-complexes and an
uncountable family of finite-type CW-complexes.

\emph{Outline of the paper:}
In Section~\ref{sectionAH} we will translate the cell attachment
problem to a purely algebraic problem using Adams-Hilton models.
We then prove our main algebraic results in Section~\ref{sectionDGA}.
In Section~\ref{sectionTop} we translate our algebraic results to
obtain Theorem~\ref{thmA}.
In Section~\ref{sectionHurewicz} we prove results about Hurewicz
images and homotopy-type to obtain Theorem~\ref{thmB}.
Finally in Section~\ref{sectionEg} we apply our results to study some
examples.

\emph{Acknowledgments:}
The results of this paper were part of my thesis at the University of
Toronto~\cite{bubenikThesis}.
I would like to thank my advisor, Paul Selick, for his encouragement
and support, and Jean-Michel Lemaire for his helpful suggestions.

\section{Adams-Hilton models}  \label{sectionAH}

Let $R = \Fp$ with $p>3$ or let $R$ be a subring of $\Q$ containing
$\frac{1}{6}$. 

A simply-connected space $X$ has an Adams-Hilton
model~\cite{adamsHilton} which we denote $\AH{X}$.
$\AH{X}$ is a connected differential graded algebra (dga) which comes
with a chain map $\AH{X} \to C_*(\lX;R)$ which induces
an isomorphism of algebras $H(\AH{X}) \isomto H_*(\lX;R)$.

Given a cell attachment $f:W \to X$ where $W=\bigvee_{j \in J}
S^{n_j}$ and $f = \vee_{j \in J} \alpha_j$, let $Y$ be the adjunction
space $Y = X \cup_f \left( \bigvee_{j \in J} e^{n_j+1} \right)$.
It is a property of Adams-Hilton models that one can choose the
following Adams-Hilton model for $Y$: 
\[
\AH{Y} = \AH{X} \amalg \freeT \langle y_j \rangle_{j \in J},
\]
where $\freeT$ denotes the tensor algebra.
The differential on $y_j$ is determined using the attaching map
$\alpha_j$.

Filter $\AH{Y}$ by the `length in $y_k$'s' filtration.
That is, let $F_{-1}\AH{Y} = 0$, let $F_0\AH{Y} = \AH{X}$, and
for $i\geq 0$, let $F_{i+1}\AH{Y} = F_i\AH{Y} + \sum_{k=0}^{i} F_k\AH{Y} \cdot
R\{y_j\}_{j \in J} \cdot F_{i-k}\AH{Y}$.
This filtration makes $\AH{Y}$ a filtered dga.

This filtration induces a first quadrant multiplicative spectral
sequence with $E^0_{p,q} = [F_p\AH{Y} / F_{p-1}\AH{Y}]_{p+q}$ which
converges from $\gr(\AH{Y})$ to $\gr(H\AH{Y})$. 

Assume that $H_*(\lX;R) \isom UL_X$ as algebras and that it is
$R$-free.
Then $(E^1, d^1) \isom (UL_X \amalg U\freeL \langle y_j \rangle_{j \in
J}, d')$, where $\freeL$ denotes the free Lie algebra and $d'$ is
determined by the induced map 
$d': R\{y_j\} \xto{d} Z\AH{X} \onto H_*(\lX;R) \isomto UL_X$.
It follows from the definition of the differential that $d'y_j =
h_X(\hat{\alpha}_j)$ where $\hat{\alpha}_j$ is the adjoint of
$\alpha_j$.
Therefore $d'y_j \in L_X$.
Thus
\[
(E^1,d^1) \isom U\eL \text{, where } 
\eL = (L_X \amalg \freeL \langle y_j \rangle, d').
\]

\section{Differential Graded Algebra Extensions} \label{sectionDGA}

Let $R = \Q$ or $\Fp$ where $p>3$ or let $R$ be a subring of $\Q$
containing $\frac{1}{6}$. 
All of our $R$-modules are graded and will be assumed to have finite type.

If $R \subset \Q$, then let $P$ be the set of invertible primes
in $R$ and let 
\[ \nP = \{p \in \Z \ | \ p \text{ is prime and } p \notin P\} \cup \{0\}.
\] 
\begin{notn} \label{notnBar}
Let $\F_0$ denote $\Q$.
If $M$ is an $R$-module then for each $p \in \nP$ we will denote $M
\tensor \Fp$ by $\bar{M}$ omitting $p$ from the notation.
Similarly if $d$ is a differential we will denote $d \tensor \Fp$ by
$\bar{d}$.
\end{notn}

Let $({A},{d})$ be a connected finite-type differential
graded algebra (dga) over $R$ which is $R$-free.
Let $Z{A}$ denote the subalgebra of cycles of ${A}$.
Let $V_1$ be a connected finite-type free $R$-module and let $d: V_1
\to Z{A}$.
Then there is a canonical dga extension $\bB = (A \amalg \freeT V_1,
d)$.

Assume that for some Lie algebra $L_0$ which is a free $R$-module,
$H({A},{d}) \isom UL_0$ as algebras.  
There is an induced map $d': V_1 \xto{d} Z{A} \to H{A}
\isomto UL_0$.  
Assume that $L_0$ can be chosen such that $d' V_1 \subset L_0$.

Taking $d'L_0 = 0$, there is a canonical differential graded Lie
algebra (dgL) 
\[ \eL = (L_0 \amalg \freeL V_1, d').
\]
Then $\eL$ is a bigraded dgL where the usual grading is called
dimension and a second grading, called degree, is given by taking
$L_0$ and $V_1$ to be in degrees $0$ and $1$ respectively.
Then the differential $d$ has bidegree $(-1,-1)$. 
\begin{notn}
Subscripts of bigraded objects will denote degree, eg. $M_0$ is the
component of $M$ in degree $0$. 
\end{notn}

The following lemma is a well-known fact, and the subsequent
lemma is part of lemmas from~\cite{anickCat2}.
We remind the reader that all of our $R$-modules have finite type.

\begin{lemma} \label{lemmaIsom}
Let $R \subset \Q$.
A homomorphism $f:M \to N$ is an isomorphism if and only if \feachpnP, $f
\tensor \Fp$ is an isomorphism.
\end{lemma}

Let $L$ be a connected bigraded dgL.
The inclusion $L \incl UL$ induces a natural map 
\begin{equation} \label{eqnPsi} 
\psi: UHL \to HUL. 
\end{equation}

\begin{lemma}[{\cite[Lemmas 4.1 and 4.3]{anickCat2}}] \label{lemmaPsi} 
Let $R = \Fp$ with $p>3$ or let $R \subset \Q$ containing $\frac{1}{6}$.
Suppose that $HUL$ is $R$-free in degrees $0$ and $1$.
Then $HL$ is $R$-free in degrees $0$ and $1$ and the map $\psi$ in
\eqref{eqnPsi} is an isomorphism in degrees $0$ and $1$.
\end{lemma}

$\bB$ is a filtered dga under the increasing filtration given by 
$F_{-1}\bB = 0$, $F_0\bB = {A}$, and for $i\geq 0$, $F_{i+1}\bB =
\sum_{j=0}^i F_{j}\bB \cdot V_1 \cdot F_{i-j}\bB$.
Letting $E^0_{p,q}(\bB) = [F_p(\bB)/F_{p-1}(\bB)]_{p+q}$ gives a first
quadrant spectral sequence of algebras: 
\[  
(E^0(\bB),d^0) = \gr(\bB) \implies E^{\infty} = \gr(H\bB).
\]

It is easy to check that $(E^1,d^1) \isom U\eL$ and hence $E^2 \isom HU\eL$. 
The following theorem follows from the main result of Anick's
thesis \cite[Theorem 3.7]{anickThesis}.
Anick's theorem holds under either of two hypotheses. 
We will use only one of these.

Recall that the Hilbert series of an $\F$-module is given by the power
series $A(z) = \Sigma_{n=0}^{\infty} (\Rank_{\F} A_n) z^n$.
Assuming that $A_0 \neq 0$, the notation $(A(z))^{-1}$ denotes the
power series $1/(A(z))$. 

\begin{thm} \label{thmMultF}
Let $R = \F$.
If the two-sided ideal $(d' V_1) \subset UL_0$ is a free $UL_0$-module
then the above spectral sequence collapses at the $E^2$ term.
That is, $\gr(H\bB) \isom HU\eL$ as algebras.
Furthermore the multiplication map
\[ \nu: \freeT(\psi(H\eL)_1) \tensor (HU\eL)_0 \to HU\eL
\]
is an isomorphism and $(HU\eL)_0 \isom UL_0/(d' V_1)$.
In addition,
\begin{multline} \label{eqnAnickFormula} 
H\bB(z)^{-1} = HU\eL(z)^{-1} 
= (1+z)(HU\eL)_0(z)^{-1} - z(UL_0)(z)^{-1} - V_1(z).
\end{multline}
\end{thm}

\begin{proof}
\cite[Theorem 3.7]{anickThesis} shows that the spectral sequence collapses as
claimed and that the multiplication map $\freeT
W \tensor (HU\eL)_0 \to HU\eL$ is an isomorphism where $W$
is a basis for $(HU\eL)_1$ as a right $(HU\eL)_0$-module.
By Lemma~\ref{lemmaPsi} and the Poincar\'{e}-Birkhoff-Witt Theorem the
homomorphism $\psi(H\eL)_1 \tensor (HU\eL)_0 \to (HU\eL)_1$ induced by
multiplication in $HU\eL$ is an isomorphism.
So we can let $W = \psi(H\eL)_1$.

The remainder of the theorem follows directly from
\cite[Theorem 3.7]{anickThesis}.
\end{proof}

\begin{cor} \label{corRfree}
If $R \subset \Q$ and \feachpnP, the two-sided ideal
$(\bar{d}\bar{V}_1) \subset U\bar{L}_0$ is a free $U\bar{L}_0$-module,
then $H\bB$ is $R$-free if and only if $HU\eL$ is $R$-free if and only
if $L_0/[d' V_1]$ is $R$-free.  
\end{cor}

\begin{proof}
First $(HU\eL)_0 \isom UL_0/(d' V_1) \isom U(L_0/[d' V_1])$.
So $(HU\eL)_0$ is $R$-free if and only if $L_0/[d' V_1]$ is $R$-free.
Since $UL_0$ and $V_1$ are $R$-free, the corollary follows from~\eqref{eqnAnickFormula}.
\end{proof}

We now prove a version of Theorem~\ref{thmMultF} for subrings of $\Q$.

\begin{thm} \label{thmMultR}
Let $R \subset \Q$.
If $L_0/[d' V_1]$ is $R$-free and \feachpnP, the two-sided ideal
$(\bar{d}\bar{V}_1) \subset U\bar{L}_0$ is a free $U\bar{L}_0$-module,
then $H\bB$ is $R$-free and the multiplication map
\[ \nu: \freeT(\psi(H\eL)_1) \tensor (HU\eL)_0 \to HU\eL
\]
is an isomorphism.
Also $\gr(H\bB) \isom HU\eL$ as algebras and
$(HU\eL)_0 \isom UL_0/(d' V_1)$. 
\end{thm}

\begin{proof}
Since $L_0/[d' V_1]$ is $R$-free, by Corollary~\ref{corRfree} so are
$HU\eL$ and $H\bB$.
It follows from the Universal Coefficient Theorem that $\fpnP$, $H\bB \tensor
\Fp \isom H(\bB \tensor \Fp)$ and $HU\eL \tensor \Fp \isom HU(\eL
\tensor \Fp)$. 
In particular $\fpnP$, $(HU\eL)_0 \tensor \Fp \isom (HU(\eL\tensor \Fp))_0$. 
Using Lemma \ref{lemmaPsi},
\[ \psi(H\eL)_1 \tensor \F_p \isom (H\eL)_1 \tensor \Fp \isom
H(\eL\tensor \Fp)_1 \isom \psi H(\eL\tensor\Fp)_1.
\]
Thus $\fpnP$, 
\[ \nu \tensor \Fp: \freeT(\psi(H\eL)_1 \tensor \Fp) \tensor (HU\eL)_0
\tensor \Fp \to H\bB \tensor \Fp
\]
corresponds under these isomorphisms to the multiplication map
\[ \freeT(\psi(H(\eL\tensor \Fp))_1) \tensor (HU(\eL \tensor \Fp))_0 \to
H(\bB \tensor \Fp). 
\]
But this is an isomorphism by Theorem~\ref{thmMultF}.
Therefore $\nu$ is an isomorphism by Lemma~\ref{lemmaIsom}.

The last two isomorphisms also follow from Theorem~\ref{thmMultF}.
\end{proof}

The next lemma will prove that if the Lie ideal $[d' V_1] \subset
L_0$ is a free Lie algebra then the hypothesis in Anick's Theorem
(Theorem~\ref{thmMultF}) holds. That is, $(d' V_1)$ is a free $UL_0$-module.

\begin{lemma} \label{lemmaI}
Given a dgL $L$ over a field $\F$, denote $UL$ by $A$.  
Let $J$ be a Lie ideal of $L$ which is a free Lie algebra, $\freeL
W$.
Take $I$ to be the two-sided ideal of $A$ generated by $J$.
Then the multiplication maps $A \tensor W \to I$ and $W \tensor A \to
I$ are isomorphisms of left and right $A$-modules respectively.
\end{lemma}

\begin{proof}
From the short exact sequence of Lie algebras 
\[ 0 \rightarrow J \rightarrow L \rightarrow L/J \rightarrow 0 \]
we get the short exact of sequence of Hopf algebras 
\[ \F \to U(J) \rightarrow U(L) \rightarrow U(L/J) \to \F \]
and so $UL \isom UJ \tensor U(L/J)$ as $\F$-modules. 
Since $J$ is a free Lie algebra $\mathbb{L}W$, $UJ \cong \freeT W$. It
is also a basic fact that $U(L/J) \isom UL/I$.
Hence we have that 
\begin{equation} \label{eqnAvs1} A \isom TW \tensor A/I 
\end{equation}
as $\F$-modules.
Furthermore
\begin{equation} \label{eqnAvs2} A \isom I \oplus A/I 
\end{equation}
as $\F$-modules.

Let $M(z)$ denote the Hilbert series for the $\F$-module $M$, and to
simplify the notation let $B = A/I$.
Then from equations \eqref{eqnAvs1} and \eqref{eqnAvs2} we have the
following (using $(TW)(z) = 1/(1-W(z))$).
\[ B(z) = A(z) (1 - W(z)), \quad
I(z) = A(z) - B(z).
\]
Combining these we have $I(z) = A(z) W(z)$. That is, $I \isom A \tensor W$
as $\F$-modules.

Let $\mu: A \tensor W \to I$ be the multiplication map. To show that
it is an isomorphism it remains to show that it is either injective
or surjective.

We claim that $\mu$ is surjective.
Since $I$ is the ideal in $A$ generated by $W$, any $x \in I$ can be
written as 
\begin{equation} \label{eqnAwb}
x= \Sigma_i a_i w_i b_{i_1} \cdots b_{i_{m_i}} \text{, where }
a_i \in A, \ w_i \in W \text{ and } b_{i_k} \in L.
\end{equation}
Each such expression gives a sequence of numbers $\{m_i\}$.
Let $M(x) = \min \left[ \max_i(m_i) \right]$, where the minimum is
taken over all possible ways of writing $x$ as in~\eqref{eqnAwb}.
We claim that $M(x) = 0$.

Assume that $M(x) = t>0$.
Then $x = x' + \Sigma_i a_i w_i b_{i_1} \cdots b_{i_t}$, where $M(x')
<t$. 
Now $w_i b_{i_1} = [w_i, b_{i_1}] \pm b_{i_1} w_i$. 
Furthermore since $J$ is a Lie ideal $[w_i, b_{i_1}] \in J \isom 
\freeL W$, so  
\[ [w_i, b_{i_1}] = \Sigma_j c_j [[w_{j_1}, \ldots , w_{j_{n_j}}] 
= \Sigma_k d_k w_{k_1} \cdots w_{k_{N_k}} 
= \Sigma_l a_l w_l,
\]
where $a_l \in A$ and $w_l \in W$.
So $x = x' + \Sigma_i \Sigma_{l_i} a_i a_{l_i} w_{l_i} b_{i_2} \cdots
b_{i_t}$.
But this is of the form in~\eqref{eqnAwb} and shows that $M(x)<t$
which is a contradiction. 

Therefore for each $x \in I$, $M(x)=0$ and we can write $x = \Sigma_i a_i
w_i$ where $a_i \in A$ and $w_i \in W$.
Then $x \in \im(\mu)$ and hence $\mu$ is an isomorphism.

Since $A$ is associative, $\mu$ is a map of left $A$-modules.

The second isomorphism follows similarly.
\end{proof}

We are now almost ready to prove our main algebraic results.
Recall that $\bB = ({A} \amalg \freeT V_1, d)$ where $d{A} \subset
{A}$ and $dV_1 \subset Z{A}$.
Also $H({A},{d}) \isom UL_0$ as algebras and if
$d':V_1 \to UL_0$ is the induced map then $d'(V_1) \subset L_0$.
Let $\eL = (L_0 \amalg \freeL V_1,d')$ with $d' L_0 = 0$.

We introduce the following terminology.
\begin{defn} \label{defnFreeDGAextension}
If $R$ is a field say $\bB$ is \emph{free} dga extension if the Lie
ideal $[d'V_1] \subset L_0$ is a free Lie algebra.
If $R \subset \Q$ say $\bB$ is \emph{free} dga extension if
$L_0/[d'V_1]$ is $R$-free and, using Notation~\ref{notnBar}, 
for every $p \in \nP$, the Lie ideal 
$[\bar{d'} \bar{V}_1] \subset \bar{L}_0$ is a free Lie algebra.
\end{defn}
\begin{defn} \label{defnSemiInertDGAextension}
In either of the cases of the previous definition we say $\bB$ is a
\emph{semi-inert} dga extension if in addition there is a free
$R$-module $K$ such that $(H\eL)_0 \sdp \freeL (H\eL)_1 \isom (H\eL)_0
\amalg \freeL K$. 
At the end of this section we will give two simpler equivalent
conditions (see Lemma~\ref{lemmaSemiInertDGA}).
\end{defn}

Note that it follows from \cite[Theorem 3.3]{halperinLemaireInert} and
\cite[Theorem 1]{felixThomasAttach} that the semi-inert condition is a
generalization of the inert condition.

Recall that there is a map $\underline{\psi}: UH\eL \to HU\eL$.
$\bB$ is a filtered dga under the increasing filtration given by
$F_{-1}\bB = 0$, $F_0\bB = {A}$, and for $i\geq 0$, $F_{i+1} =
\sum_{k=0}^i F_k\bB \cdot V_1 \cdot F_{i-k}\bB$.
There is an induced filtration on $H\bB$.
We prove one last lemma.

\begin{lemma} \label{lemmaQuotientf}
There exists a quotient map
\begin{equation} \label{eqnQuotientf} 
f: F_1 H\bB \to (HU\eL)_1.
\end{equation}
Given $\bar{w} \in (H\eL)_1$ there exists a cycle
$w \in F_1 \bB$ such that $f([w])=\bar{w}$.
\end{lemma}

\begin{proof}
By Theorem~\ref{thmMultF} or Theorem~\ref{thmMultR}, 
$(\gr(H\bB))_1 \isom (HU\eL)_1$. 
So there is a quotient map 
\[ f: F_1 H\bB \onto (\gr(H\bB))_1 \isomto (HU\eL)_1.
\]
By Lemma~\ref{lemmaPsi} $(H\eL)_1 \isom (\underline{\psi}H\eL)_1
\subset (HU\eL)_1$.
So for $\bar{w}$ one can choose a representative cycle $w \in
Z F_1 \bB$ such that $f([w]) = \bar{w}$.
\end{proof}

Recall that $R = \Fp$ with $p>3$ or $R \subset \Q$ containing
$\frac{1}{6}$.
Also all of our $R$-modules are connected, $R$-free, and have finite
type.
Let $({A}, {d})$ be a dga and let $V_1$ be a $R$-module with a map
$d:V_1 \to {A}$. 
Assume that there exists a Lie algebra $L_0$ such that
$H({A},{d}) \isom UL_0$ as algebras and $d' V_1 \subset
L_0$ where $d'$ is the induced map.

\begin{thm} \label{thmDGAextension}
Let $\bB = ({A} \amalg \freeT V_1, d)$.
Assume that $\bB$ is a \emph{free} dga extension in the sense of
Definition~\ref{defnFreeDGAextension}. 
Let $\eL = (L_0 \amalg \freeL V_1, d')$. \\ 
(a) Then as algebras
\[ \gr(H\bB) \isom U ((H\eL)_0 \sdp \freeL (H\eL)_1)
\]
with $(H\eL)_0 \isom L_0/[d' V_1]$ as Lie algebras.  
If $R \subset \Q$ then additionally 
$(H\eL)_0 \sdp \freeL (H\eL)_1 \isom \underline{\psi} H\eL$ 
as Lie algebras. \\
(b) Furthermore if $\bB$ is \emph{semi-inert} (that is, there is an
$R$-module $K$ such that $(H\eL)_0 \sdp \freeL (H\eL)_1 \isom (H\eL)_0
\amalg \freeL K$) then as algebras 
\[ H\bB \isom U ( (H\eL)_0 \amalg  \freeL K' )
\]
for some $K' \subset F_1 H\bB$ such that $f: K' \isomto K$, where
$f$ is the quotient map in Lemma~\ref{lemmaQuotientf}.
\end{thm}

\begin{proof}
(a) If $R=\Fp$ then by Lemma~\ref{lemmaI}, $(d' V_1) \subset UL_0$ is
a free $UL_0$-module. 
If $R\subset \Q$ then by Lemma~\ref{lemmaI}, for each $p \in \nP$,
$(\bar{d}' \bar{V}_1) \subset U\bar{L}_0$ is a free $U\bar{L}_0$-module. 
So we can apply either Theorem~\ref{thmMultF} or
Theorem~\ref{thmMultR} to show that $\gr(H\bB) \isom HU\eL$ as
algebras and that the multiplication map
\[ \nu: \freeT(\underline{\psi}(H\eL)_1) \tensor (HU\eL)_0 \to HU\eL
\]
is an isomorphism.

By Lemma~\ref{lemmaPsi} $(HU\eL)_0 \isom U(H\eL)_0$ and
$\underline{\psi}(H\eL)_1 \isom (H\eL)_1$.
By the definition of homology $(H\eL)_0 \isom L_0/[d' V_1]$.

Let $N = \underline{\psi}(H\eL)$.
Then $N_0$ acts on $N_1$ so we can define
\[
L' = N_0 \sdp \freeL N_1.
\]
Note that $L'_0 = N_0$ and $L'_1 = N_1$.
There is a Lie algebra map $u: L' \to N$ and
an induced algebra map $\tilde{u}: UL' \xrightarrow{Uu}
UN \to HU\eL$. 

Recall that as $R$-modules, $L' \isom N_0 \times \freeL N_1$.
The Poincar\'{e}-Birkhoff-Witt Theorem shows that the
multiplication map
\[ \phi: \freeT N_1 \tensor (HU\eL)_0 \isomto U\freeL N_1 \tensor U N_0 \to
UL'
\]
is an isomorphism.
Since $\tilde{u}$ is an algebra map, $\nu = \tilde{u} \circ \phi$.
Thus $\tilde{u}$ is an isomorphism.
Therefore $HU\eL \isom UL'$ as algebras and hence $\gr(H\bB) \isom
UL'$.
If $R=\Fp$ then this finishes (a).

If $R \subset \Q$ then we will show that $u: L' \to N$ is an isomorphism. 
Let $\iota:N \incl HU\eL$ be the inclusion.
Since the composition $L' \xrightarrow{u} N
\xrightarrow{\iota} HU\eL \isomto UL'$ is 
the canonical inclusion $L' \incl UL'$, $u$ is injective.
The inclusion $L' \to UL'$ splits as $R$-modules;
so as $R$-modules $N \isom L' \oplus N/L'$.
Since $L'$ and $N$ are $R$-free, so is $N/L'$.

Recall that the composition $\iota \circ u$ induces the isomorphism
$\tilde{u}: UL' \xrightarrow{Uu} UN \to HU\eL$.
Tensor these maps with $\Q$ to get the commutative diagram
\begin{equation} \label{eqnCdq}
\xymatrix{
UL'\tensor\Q \ar[r]^{Uu\tensor\Q} \ar[dr]_{\isom} & UN \tensor \Q
\ar[d] \\
& HU\eL \tensor \Q }.
\end{equation}
It is a classical result that the natural map
\begin{equation} \label{eqnPsiq} 
\underline{\psi}_{\Q}: UH(\eL \tensor \Q) \xrightarrow{\isom} HU(\eL
\tensor \Q)
\end{equation} 
is an isomorphism. 
Notice that
\[ N \tensor \Q = (\underline{\psi}H\eL) \tensor \Q
\isom \underline{\psi}_{\Q} H(\eL \tensor \Q) \isom H(\eL \tensor \Q)
\]
and $HU\eL \tensor \Q \isom HU(\eL \tensor \Q)$.
Under these isomorphisms the vertical map in~\eqref{eqnCdq}
corresponds to the isomorphism in~\eqref{eqnPsiq}.

Therefore $Uu \tensor \Q$ is an isomorphism and hence $u\tensor\Q$ is
surjective.
As a result $\coker u = N/L'$ is a torsion $R$-module.
But we have already shown that $N/L'$ is $R$-free.
Thus $N/L' = 0$ and $N \isom
L'$. 
Hence $HU\eL \isom UN$.

\noindent 
(b) Recall that $N_0$ acts on $N_1  = (\underline{\psi}H\eL)_1
\isom (H\eL)_1$ via the adjoint action.
Assume that $\bB$ is semi-inert. That is, there exists $\{\bar{w}_i\} \subset N_1$
such that $L' \isom N_0 \amalg \freeL K$, where $K = R\{\bar{w}_i\}$. 
Recall from (a) that $HU\eL \isom \gr(H\bB)$ and that the inclusions
$N_0 \subset HU\eL$ and $\bar{w}_i \in HU\eL$ induce a Lie algebra map 
\[ u: L' \to \gr(H\bB).
\]

By Lemma~\ref{lemmaQuotientf}, there exists $w_i \in F_1 \bB$ such that
$f([w_i]) = \bar{w}_i$ where $f$ is the map in~\eqref{eqnQuotientf}.
Let $K' = R\{[w_i]\} \subset F_1 H\bB$, and let $L'' = N_0 \amalg
\freeL K'$. 
Then $f: K' \xrightarrow{\isom} K$ and $f$ induces an isomorphism
$L'' \xrightarrow{\isom} L'$.

By part (a), $N_0 \subset (\gr H\bB)_0$. 
Since $F_{-1} H\bB = 0$, $(\gr H\bB)_0 \incl H\bB$, so $N_0 \incl H\bB$.
Since $N_0$, $K' \incl H\bB$, there are induced maps 
\[
\xymatrix{
L'' \ar[r]^-{\eta} \ar[d] & H\bB \\
UL'' \ar[ur]_{\theta}
},
\]
where $\eta$ is a Lie algebra map and $\theta$ is an algebra map.

Grade $L''$ by letting $N_0$ be in degree $0$ and $K'$ be in
degree $1$.
This also filters $L''$.
Then $\eta$ is a map of filtered objects.

From this we get the following commutative diagram
\[
\xymatrix{
\gr(L'') \ar[d] \ar[r]^-{\gr(\eta)} & \gr(H\bB) \\
\gr(UL'') \ar[ur]^{\gr(\theta)} \ar[d]_{\isom} \\
U\gr(L'') \ar[uur]_{\rho}
}.
\]

Now $\gr(L'') \isom L'' \isom L'$ and $\gr(\eta)$ corresponds to $u$ under
this isomorphism. 
So $\rho$ corresponds to $\tilde{u}$ which is an isomorphism.
Thus $\gr(\theta)$ is an isomorphism, and hence $\theta$ is an
isomorphism. 
Therefore $H\bB \isom UL''$ which finishes the proof.
\end{proof}

As promised we now give two simpler equivalent conditions for
semi-inertness.

\begin{lemma} \label{lemmaSemiInertDGA}
Let $\bB$ be a free dga extension (in the sense of
Definition~\ref{defnFreeDGAextension}). 
Then the following conditions are equivalent: \\
(a) $(H\eL)_0 \sdp \freeL (H\eL)_1 \isom (H\eL)_0 \amalg \freeL K$ 
for some free $R$-module $K \subset (H\eL)_1$, \\
(b) $(H\eL)_1$ is a free $(H\eL)_0$-module, and \\
(c) $\gr_1(H\bB)$ is a free $\gr_0(H\bB)$-bimodule.
\end{lemma}

\begin{proof}
(b) $\implies$ (a)
Let $K$ be a basis for $(H\eL)_1$ as a free $(H\eL)_0$-module.
Then $(H\eL)_0 \sdp \freeL (H\eL)_1 \isom (H\eL)_0 \amalg \freeL
K$. \\
(a) $\implies$ (c)
Since $\bB$ is a free dga extension, by Theorem~\ref{thmDGAextension}(a),
$\gr_*(H\bB) \isom U \left( (H\eL)_0 \sdp \freeL (H\eL)_1 \right)$.
So by (a), 
\[
\gr_*(H\bB) \isom U \left( (H\eL)_0 \amalg \freeL K \right) \isom
\gr_0(H\bB) \amalg \freeT K,
\]
for some free $R$-module $K \subset (H\eL)_1$. 
Therefore 
\[
\gr_1(H\bB) \isom \left[ \gr_0(H\bB) \amalg \freeT K \right]_1 \isom
\gr_0(H\bB) \tensor K \tensor \gr_0(H\bB).
\]
(c) $\implies$ (b)
Let $L' = (H\eL)_0 \sdp \freeL (H\eL)_1$.
Then by Theorem~\ref{thmDGAextension}(a), $\gr_*(H\bB) \isom UL'$ and
$\gr_1(H\bB) \isom (UL')_1$. 
By (c), $(UL')_1$ is a free $(UL')_0$-bimodule.
Then we claim that it follows that $N_1$ is a free $N_0$-module.
Indeed, if there is a nontrivial degree one relation in $L'$ then
there is a corresponding nontrivial degree one relation in $UL'$.
\end{proof}

\section{Cell-attachments} \label{sectionTop}

Let $R = \Fp$ with $p > 3$ or $R \subset \Q$ containing $\frac{1}{6}$.
Let $X$ be a finite-type simply-connected CW-complex such that
$H_*(\lX;R)$ is torsion-free and as algebras $H_*(\lX;R) \isom UL_X$
where $L_X$ is the Lie algebra of Hurewicz images. 
Let $W = \bigvee_{j \in J} S^{n_j}$ be a finite-type wedge of spheres and
let $f:W \to X$.
Let $Y = X \cup_f \bigl(\bigvee_{j \in J} e^{n_j+1}\bigr)$.
Using the Adams-Hilton models of Section~\ref{sectionAH}, we defined
a Lie algebra $\eL = (L_X \amalg \freeL \langle y_j \rangle, d')$.
Mirroring Definitions \ref{defnFreeDGAextension} and
\ref{defnSemiInertDGAextension}, we introduce the following terminology.

\begin{defn} \label{defnFreeCA}
If $R$ is a field call $f$ a \emph{free} cell attachment if the Lie
ideal $[L^W_X] \subset L_X$ is a free Lie algebra.
If $R \subset \Q$ call $f$ a \emph{free} cell attachment if
$L_X/[L^W_X]$ is $R$-free and for every $p \in \nP$, the Lie ideal
$[\bar{L}^W_X] \subset \bar{L}_X$ is a free Lie algebra.
\end{defn}
\begin{defn} \label{defnSemiInertCA}
In either of the cases of the previous definition say that $f$ is a
\emph{semi-inert} cell attachment if in addition one of the following
three equivalent (by Lemma~\ref{lemmaSemiInertDGA}] conditions is
satisfied: \\ 
(a) $\gr_1(H\AH{Y})$ is a free $\gr_0(H\AH{Y})$-bimodule, \\
(b) $(H\eL)_1$ is a free $(H\eL)_0$-module, or \\
(c) there is a free $R$-module $K$ such that 
\[ (H\eL)_0 \sdp \freeL (H\eL)_1 \isom (H\eL)_0 \amalg \freeL K.
\]
\end{defn}

Theorem \ref{thmDGAextension} gives most of the following topological
result directly. 

\begin{thm}  \label{thmPreciseA}
Let $Y = X \cup_f \bigl(\bigvee_{j \in J} e^{n_j+1}\bigr)$ and let
$\eL = (L_X \amalg \freeL \langle y_j \rangle, d')$.
Assume that $f$ is \emph{free}. \\
(a) Then $H_*(\lY;R)$ and $\gr(H_*(\lY;R))$ are $R$-free and as algebras 
\[ \gr(H_*(\lY; R)) \isom U( {L}^X_Y \sdp \freeL (H\eL)_1)
\]
with ${L}^X_Y \isom {L}_X/[{L}^W_X]$ as Lie algebras. \\ 
(b) Furthermore if $f$ is \emph{semi-inert} then as algebras 
\[ H_*(\lY;R) \isom U ( L^X_Y \amalg \freeL K' )
\]
for some $K' \subset F_1 H_*(\lY; R)$. 
\end{thm}

\begin{proof}
It remains to show that $(H\eL)_0 \isom L^X_Y$.
By Theorem~\ref{thmDGAextension} we have the algebra isomorphism
\[
\gr(H\AH{Y}) \isom U( (H\eL)_0 \sdp \freeL (H\eL)_1)
\]
with $(H\eL)_0 \isom L_X / [L^W_X]$.
Therefore 
\begin{equation} \label{eqnF0hul} 
F_0 H\AH{Y} \isom (\gr(H\AH{Y}))_0 \isom U(H\eL)_0 \isom U(L_X/[L^W_X]). 
\end{equation}

The inclusion $i: \AH{X} \isomto F_0 \AH{Y}$ induces a map $H(i):
H\AH{X} \to F_0 H\AH{Y}$.
Now under the isomorphism \eqref{eqnF0hul} and $UL_X \isomto
H\AH{X}$ the  map $H(i)$ corresponds to a map $UL_X \to
U(L_X/[L^W_X])$ where $U(L_X/[L^W_X]) \subset UL_Y$.
It is easy to check that this is the canonical map.
In other words $L^X_Y \isom L_X/[L^W_X]$.
Therefore $(H\eL)_0 \isom L^X_Y$.
\end{proof}

Corollary~\ref{corNice} follows from Theorem \ref{thmPreciseA}.

\begin{proof}[Proof of Corollary~\ref{corNice}]
The cell attachment $f$ is nice in the sense of Hess and
Lemaire~\cite{hessLemaireNice} if and only if $UL_X/(L^W_X)$ injects in
$H_*(\lY;R)$.
Recall the standard fact that $UL_X/(L^W_X) \isom U(L_X/[L^W_X])$.
By Theorem~\ref{thmPreciseA},
\[ U(L_X/[L^W_x]) \isom UL^X_Y \isom \gr_0 (H_*(\lY;R))
\]
which injects in $H_*(\lY;R)$.
\end{proof}

\section{Hurewicz images} \label{sectionHurewicz}

Let $R = \Fp$ with $p >3$ or $R \subset \Q$ be a subring containing
$\frac{1}{6}$. 

Recall that we have a homotopy cofibration $W \xto{f} X \to Y$ where $W =
\bigvee_{j \in J} S^{m_j}$ is a finite-type wedge of spheres, $f =
\bigvee_{j \in J} \alpha_j$, $H_*(\lX;R)$ is torsion-free, and as
algebras $H_*(\lX;R) \isom UL_X$. 
Let $\hat{\alpha}$ denote the adjoint of $\alpha$.

Assume that $f$ is free. That is, the Lie ideal $[L^W_X]$ is a free
Lie algebra. 

Recall that $h_X: \pi_*(\lX) \tensor R \to L_X \subset H_*(\lX;R)$ is
a Lie algebra map.
In order to identify Hurewicz images in $H_*(\lY;R)$ we will need to
be able to construct maps from information about the loop space
homology.
In particular, we will need to assume that there exists a Lie algebra
map $\sigma_X: L_X \to \pi_*(\lX)\tensor R$ such that $h_X \circ
\sigma_X = id_{L_X}$.  
This map exists if $R = \Q$ or if $X$ is a wedge of spheres.
We will give a sufficient condition for the existence of this map
later in this section.

If $R$ is a subring of $\Q$ with invertible primes $P \supset 
\{2,3\}$, then we may need to exclude those primes $p$ for which an
attaching map $\alpha_j \in \pi_*(X)$ includes a term with $p$-torsion. 
Following Anick~\cite{anickCat2} we define the set of \emph{implicit
primes} below. 

\begin{defn} \label{defnIp}
By the Milnor-Moore theorem~\cite{milnorMooreHopfAlgebras}, $h_X,
\sigma_X$ are rational isomorphisms, so $\im(\sigma_X \circ h_X - id)$
is a torsion element of $\pi_*(\lX) \tensor R$.
Let $\gamma_j = \sigma_X h_X(\hat{\alpha_j}) - \hat{\alpha}_j$ where
$\hat{\alpha}_j: S^{n_j-1} \to \lX$ is the adjoint of $\alpha_j$.
Then $t_j \gamma_j = 0 \in \pi_*(\lX) \tensor R$ for some $t_j>0$.
Let $t_j$ be the smallest such integer.
Define $P_Y$, the set of \emph{implicit primes} of $Y$ as follows.
A prime $p$ is in $P_Y$ if and only if $p \in P$ or 
$p | t_j$ for some $j \in J$.
\end{defn}

One can verify that the implicit primes have the following properties.

\begin{lemma} \label{lemmaIp}
(a) Let $\{x_i\}$ be a set of Lie algebra generators for $L_X$ and let
$\beta_i = \sigma_X x_i$.
If all of the attaching maps are $R$-linear combinations of
iterated Whitehead products of the maps $\beta_i$, then $P_Y = P$. \\
(b) If $P = \{2,3\}$ and $n = \dim (Y)$ then the implicit primes are
bounded by $\max (3,n/2)$. 
\end{lemma}




By replacing $R$ with $R'=\Z[{P_Y}^{-1}]$ if necessary, we may assume
that the implicit primes are invertible.
This implies that for all $j \in J$, $\sigma_X h_X
\hat{\alpha}_j = \hat{\alpha}_j$,
and hence $\sigma_X(dy_j) = \hat{\alpha_j}$.

\begin{rem}
If $R = \Fp$ then we will also need that for all $j \in J$, 
$\sigma_X h_X \hat{\alpha}_j = \hat{\alpha}_j$. 
So that we can state the cases $R=\Fp$ and $R\subset \Q$ together,
when $R=\Fp$ and we say \emph{the implicit primes are invertible} 
we mean that for all  
$j \in J$,  $\sigma_X h_X \hat{\alpha}_j = \hat{\alpha}_j$. 
\end{rem}

We now consider both cases $R = \Fp$ or $R\subset \Q$.
Recall that in Section~\ref{sectionAH} we defined the differential
graded Lie algebra 
\[ \eL = (L_X \amalg \freeL \langle y_j \rangle_{j   \in J}, d').
\]
By Theorem \ref{thmPreciseA}, $H_*(\lY;R)$ is torsion-free
and $\gr(H_*(\lY;R)) \isom U (L^X_Y \sdp \freeL (H\eL)_1)$ as algebras.
From this we want to show that $H_*(\lY;R) \isom UL_Y$ as algebras.

This situation closely resembles that of torsion-free spherical two-cones,
and we will generalize Anick's proof for that situation~\cite{anickCat2}.

The proof of the following is the same as the proof
of~\cite[Claim 4.7]{anickCat2}, so we will only sketch it here.
See~\cite[Chapter 7]{bubenikThesis} for more details.

\begin{prop} \label{propU1}
Let $W \to X \to Y$ and $\eL$ be as above with $\sigma_X h_X
\hat{\alpha}_j = \hat{\alpha_j}$, $\forall j$.
Then there exists an injection of $L^X_Y$-modules
\[ u_1: (H\eL)_1 \incl (\gr(L_Y))_1.
\]
\end{prop}

\begin{proof}[Proof sketch.]
Let $\bar{w}_i \in (H\eL)_1$ be a homology class in dimension $m+1$.
It can be represented by a cycle $\gamma_i \in \eL$ in degree $1$.
Using the Jacobi identity one can write
\[ \gamma_i = \Sigma_{k=1}^s c_k u_k \text{ where } c_k \in R \text{
  and } u_k = [ \cdots [ y_{j_k}, x_{k_1}] ,\ldots ,x_{k_{n_k}}],
\]
with $[x_{k_i}] \in L_X$.

The sphere $S^m$ has an Adams-Hilton model $(\freeT\langle a \rangle,
0)$ which can be extended to $(\freeT\langle a, b \rangle,
d)$ where $db=a$ which is an Adams-Hilton model for the disk
$D^{m+1}$.
Using properties of Adams-Hilton models, one can construct maps 
$g_k: (D^{m+1}, S^m) \to (Y,X)$ for $1 \leq k \leq s$ such that 
$\AH{g_k}(b)=c_k u_k$.

Let 
\[
g' = g_1 \vee \ldots \vee g_s: \left(\bigvee_{k=1}^s D^{m+1},
\bigvee_{k=1}^s S^m \right) \to (Y,X).
\]

Using the fact that $\sigma_X h_X \hat{\alpha}_j = \hat{\alpha_j}$, one
can show that $g'|_{\bigvee_k S^m}$ is contractible in $X$.
As a result, $g'$ can be extended to a map $g: S^{m+1} \to Y$ whose
Hurewicz image modulo lower filtration is $\bar{w}_i$.
\end{proof}

From this we will prove our final main result.
Recall that $R = \Fp$ with $p >3$ or $R \subset \Q$ containing $\frac{1}{6}$.

\begin{thm} \label{thmHurewicz}
Let $Y = X \cup_f 
\left(\bigvee e^{n_j+1}\right)$
with $f=\bigvee \alpha_j$ satisfying the hypotheses of Theorem~\ref{thmPreciseA}.
In addition assume that $h_X$ has a right inverse $\sigma_X$ and that
the implicit primes are invertible.
Then \\
\noindent (a) the canonical algebra map
\begin{equation} \label{eqnUL_Y2} 
UL_Y \to H_*(\lY;R)
\end{equation} 
is a surjection. \\
\noindent (b) If $R \subset \Q$ then \eqref{eqnUL_Y2} is an
isomorphism, $\gr(L_Y) \isom L^X_Y \sdp \freeL((H\eL)_1)$, and
localized at $R$, $\lY \in \PiS$. \\
\noindent (c) If $R \subset \Q$ and $f$ is \emph{semi-inert} then \\
(i) there exists $\hat{K} \subset F_1 L_Y$ such that $L_Y
\isom L^X_Y \amalg \freeL \hat{K}$ as Lie algebras, \\
(ii) $L_Y \isom H\eL$ as Lie algebras, and \\
(iii) $h_Y$ has a right inverse $\sigma_Y$.
\end{thm}

\begin{proof}[Proof of (a) and (b)] 
\ \\ 
\noindent
(a)
Let $g: \gr(H_*(\lY;R)) \isomto UL'$ be the algebra isomorphism given by
Theorem~\ref{thmPreciseA}(a) where $L' = L^X_Y \sdp \freeL K'$ with $K' =
(H\eL)_1$. 
Note that $L'_0 = L^X_Y$ and that $L'_1 = K'$.

We have an injection of Lie algebras 
\begin{equation} \label{eqnL0prime}
u_0: L^X_Y \incl F_0 L_Y \isomto (\gr(L_Y))_0.
\end{equation}

Since the implicit primes are invertible, we
have that for all $j$, $h_X \sigma_X \hat{\alpha}_j = \hat{\alpha}_j$.
So by Proposition~\ref{propU1} we get an injection of $L^X_Y$-modules 
$u_1: K' \incl (\gr(L_Y))_1$.
Hence for $x \in L^X_Y$ and $y \in K'$, $u_1([y,x]) = [u_1(y), u_0(x)]$.
Thus $u_{0}$ and $u_{1}$ can be extended to a Lie algebra map $u: L'
\rightarrow \gr(L_Y)$. 

The inclusion $L_Y \hookrightarrow H_*(\lY;R)$ induces a map between 
the corresponding graded objects, $\chi: \gr(L_Y) \rightarrow
\gr(H_*(\lY;R))$.

We claim that for $j=0$ and $1$, $g\circ \chi \circ u_{j}$ is the ordinary 
inclusion of $L'_{j}$ in $UL'$.
For $j=0$, under the isomorphisms $\gr_0 L_Y \isom F_0 L_Y$ and $\gr_0
H_*(\lY;R) \isom F_0 H_*(\lY;R)$, $g \chi u_0$ corresponds to the map
$L^X_Y \incl F_0 L_Y \incl F_0 H_*(\lY;R) \isomto UL^X_Y$.
For $j=1$, under the isomorphism $UL' \isom HU\eL$, $g \chi u_1$
corresponds to the inclusion $K' = (H\eL)_1 \incl (HU\eL)_1$.
It follows that $g \circ \chi \circ u$ is the standard inclusion $L'
\incl UL'$. 
Since $g \circ \chi \circ u$ is an injection, so is $u$.

The canonical map $U\gr(L_Y) \isomto \gr(UL_Y)$ is an algebra isomorphism. 
Now $u$ and $\chi$ induce the maps $Uu$ and  $\bar{\chi}$
in the following diagram.
\begin{equation} \label{cdUL_Y} 
\xymatrix{
UL' \ar[r]^-{Uu} \ar[dr]^{\isom}_{g^{-1}} & U\gr(L_Y) \ar[r]^-{\cong} 
\ar[d]^{\bar{\chi}} & \gr(UL_Y) \ar@{-->}[dl]^{\tilde{\chi}} \\
& \gr(H_*(\lY;R)) }
\end{equation}
Since we showed that $g \chi u$ is the ordinary inclusion $L' \incl
UL'$ the diagram commutes.
Since $g^{-1}$ is surjective, the induced map $\tilde{\chi}$ is surjective.
Since $\tilde{\chi}$ is the associated graded map to the canonical map
$UL_Y \to H_*(\lY;R)$ and the filtrations are bicomplete, the associated
ungraded map is also surjective. 
So the canonical map $UL_Y \to H_*(\lY;R)$ is surjective which
finishes the proof of (i).
\\

\noindent
(b)
In the case where $R \subset \Q$ we can tensor with $\Q$ and make use
of results from rational homotopy theory. 

Recall that $\gr(H_*(\lY;R)) \isom UL'$ and that we constructed a Lie
algebra map $u: L' \to \gr(L_Y)$ and showed that it is an injection.
We claim that for $R \subset \Q$, $u$ is an isomorphism. 
$H_*(\lY;R)$ and $\gr(H_*(\lY;R))$ have the same Hilbert series.
Also since $H_*(\lY; R)$ is torsion-free, it has the same Hilbert
series as $H_*(\lY; \Q)$.
Let $S$ be the image of $h_Y \tensor \Q$. 
Then $S$, $L_Y$ and $\gr(L_Y)$ have the same Hilbert series.
By the Milnor-Moore Theorem~\cite{milnorMooreHopfAlgebras}, $H_*(\lY;
\Q) \isom US$. 
So by the Poincar\'{e}-Birkhoff-Witt Theorem, $S$ has the same Hilbert
series as $L'$, and hence $\gr(L_Y) \isom L'$ as $R$-modules. 
Since $u: L' \to \gr(L_Y)$ is an injection it follows that it is an
isomorphism.  

Using diagram~\eqref{cdUL_Y} we get that~\eqref{eqnUL_Y2} is an
isomorphism. 

The map $u$ gives the desired Lie algebra isomorphism $\gr(L_Y) \isom
L' = L^X_Y \sdp \freeL (H\eL)_1$.

By the Hilton-Serre-Baues 
Theorem~\cite[Lemma V.3.10]{bauesCommutator},\cite[Lemma 3.1]{anickCat2},
that \eqref{eqnUL_Y2} is an isomorphism is equivalent to the
statement that localized at $R$, $\lY \in \PiS$.
\end{proof}

Before we prove (c) we strengthen the result in (b) in the semi-inert case. 

\begin{lemma} \label{lemmaStrengthen}
Let $R \subset \Q$. 
Assume $Y$ is a space satisfying the hypotheses of
Theorem~\ref{thmHurewicz}.
If furthermore $f$ is semi-inert then there exists $\hat{K} \subset F_1
L_Y$ such that $L_Y \isom L^X_Y \amalg \freeL \hat{K}$ as Lie
algebras.
\end{lemma}

\begin{proof}
Assume that $f$ is semi-inert.
Recall the situation from Theorem~\ref{thmPreciseA}(b).
We have that $\gr(H_*(\lY;R)) \isom UL'$ where $L' \isom L^X_Y
\amalg \freeL K'$ where $K' \isom R\{\bar{w}_i\} \subset
(\gr(H_*(\lY;R)))_1$.
For each $\bar{w}_i$ let $[w_i]$ be an inverse image under the
quotient map $H_*(\lY;R) \to \gr(H_*(\lY;R))$.
Let $K'' = \{[w_i]\}$.
By Theorem~\ref{thmPreciseA}(b) as algebras $H_*(\lY;R) \isom UL''$ where
$L'' = L^X_Y \amalg \freeL K''$ (see Theorem~\ref{thmDGAextension}(b)).
Since $K'' \xrightarrow{\isom} K'$, there is an induced Lie algebra
isomorphism $L'' \xrightarrow{\isom} L'$.
So $\gr(H_*(\lY;R)) \isom H_*(\lY;R)$ as algebras. 

Recall that $K' = (H\eL)_1$, so by Proposition~\ref{propU1} there
exists $\hat{K} \subset F_1 L_Y$ such that $f: \hat{K} \isomto K'$,
where $f$ is the quotient map from~\eqref{eqnQuotientf}.
Let 
\[
\hat{L} = L^X_Y \amalg \freeL \hat{K}.
\]
So $f: \hat{K} \isomto K'$, induces a Lie algebra isomorphism
$\hat{L} \isomto L'$.
This in turn induces the algebra isomorphism $U\hat{L} \isomto UL'
\isomto H_*(\lY;R)$. 

Since $U\hat{L} \isom H_*(\lY;R)$ as algebras there is an injection
$L_Y \incl U\hat{L}$.
Also, since $\hat{K} \subset L_Y$ there is a canonical Lie algebra map
$v: \hat{L} \to L_Y$.
These fit into the following commutative diagram.
\[
\xymatrix{
\hat{L} \ar[rr] \ar[dr]_v & & U\hat{L} \\
& L_Y \ar[ur]
}
\]
It follows that $v$ is an injection.

We claim that $v$ is an isomorphism. 
Since $H_*(\lY; R)$ is torsion-free, it has the same Hilbert
series as $H_*(\lY; \Q)$.
Let $S$ be the image of $h_Y \tensor \Q$. 
Then $S$ and $L_Y$ have the same Hilbert series.
By the Milnor-Moore Theorem~\cite{milnorMooreHopfAlgebras}, $H_*(\lY;
\Q) \isom US$. 
So by the Poincar\'{e}-Birkhoff-Witt Theorem,
$S$ has the same Hilbert series as $\hat{L}$, and hence $L_Y \isom
\hat{L}$ as $R$-modules.  
Since $v: \hat{L} \to L_Y$ is an injection it follows that it is an
isomorphism.  

Therefore $L_Y \isom L^X_Y \amalg \freeL \hat{K}$ as Lie algebras,
with $\hat{K} \subset F_1 L_Y$ and
$H_*(\lY;R) \isom UL_Y$ as algebras. 
\end{proof}

\begin{proof}[Proof of Theorem~\ref{thmHurewicz}(c)] 
Lemma~\ref{lemmaStrengthen} proves (i) and puts us in a position to
prove (ii) and (iii).

(ii) 
Recall that $(H\eL)_1 \isom \hat{K}$ and  $L' = (H\eL)_0 \amalg
\freeL ( (H\eL)_1 )$.  
Thus there is a canonical map $\nu: L' \to H\eL$.
Now the composition
\[ L' \xto{\nu} H\eL \xto{\phi} HU\eL \isomto UL'
\]
is just the ordinary inclusion of L' into its universal enveloping
algebra.
Therefore $\nu$ is an injection.
Tensoring with $\Q$ we get that $L'$ and $H\eL$ have the same Hilbert
series.
It follows that $\nu$ is an isomorphism. 

(iii) 
We will construct a map $\sigma_Y$ right inverse to $h_Y$.  

Let $i$ denote the inclusion $X \incl Y$.
Consider the composite map 
\[ F: [L^W_X] \incl L_X \xrightarrow{\sigma_X} \pi_*(\lX) \tensor R
\xrightarrow{(\Omega i)_{\#}} \pi_*(\lY) \tensor R.
\]
We claim that $F=0$. Since $F$ is a Lie algebra map it is sufficient
to show that it is zero on the Lie algebra generators of $L^W_X$.
That is, show
\[ (\Omega i)_{\#} \sigma_X (R\{h_X(\hat{\alpha}_j)\}) = 0.
\]
Since there are no implicit primes $\sigma_X h_X \hat{\alpha}_j =
\hat{\alpha}_j$. 
By the construction of $Y$, $\Omega i \circ \hat{\alpha}_j \simeq 0$.
So $F=0$ as claimed.

Therefore there is an induced map $G: L^X_Y \isomto L_X/[L^W_X] \to
\pi_*(\lY) \tensor R$. 
That $h_Y \circ G$ is the inclusion map can be seen from the following
commutative diagram. 
\[
\xymatrix{
L_X \ar@/^2.35ex/@{-->}[r]^-{\sigma_X} \ar@{->>}[d] & \pi_*(\lX)
\tensor R \ar[dd]^{(\Omega i)_{\#}} \ar[l]^-{h_X} \\ 
L^X_Y \ar[dr]^-{G} \ar[d] \\
L_Y & \pi_*(\lY) \tensor R \ar[l]_-{h_Y}
}
\]

Now construct $\sigma_Y: L_Y \to \pi_*(\lY) \tensor R$ as follows.
We have shown that $L_Y \isom L^X_Y \amalg \freeL \hat{K}$ for some $\hat{K}
\subset F_1 L_Y$.  
Since 
$h_Y: \pi_*(\lY) \tensor R \onto L_Y$, choose preimages $\check{K}
\subset \pi_*(\lY) \tensor R$ such that $h_Y: \check{K} \isomto \hat{K}$.
Let $\sigma_Y \!\! \mid_{L^X_Y} = G$ and let 
$\sigma_Y \hat{K} = \check{K}$ be right inverse to $h_Y$.
Now extend $\sigma_Y$ canonically to a Lie algebra map on $L_Y$.

We finally claim that $h_Y \sigma_Y = id_{L_Y}$.
Since $h_Y \sigma_Y$ is a Lie algebra map it suffices to check that it
is the identity for the generators. 
\[ h_Y \sigma_Y L^X_Y = h_Y G L^X_Y = L^X_Y, \quad
h_Y \sigma_Y \hat{K} = h_Y \check{K} = \hat{K}
\]
Therefore $\sigma_Y$ is the desired Lie algebra map right inverse to $h_Y$.
\end{proof}

\begin{rem}
Anick conjectured Theorem~\ref{thmHurewicz}(c)(ii) without the
semi-inert condition in the special case where $X$ is a wedge of
spheres~\cite[Conj. 4.9]{anickCat2}.
\end{rem}

\begin{cor}[Corollary~\ref{corTFMM}]
If $R \subset \Q$ then the canonical algebra map 
\[ U(F(\pi_*(\lY) \tensor R)) \to H_*(\lY;R)
\]
is an isomorphism.
\end{cor}

\begin{proof}
Recall that $FM = M / \Torsion(M)$.
By definition, the Lie algebra map $h_Y: \pi_*(\lY) \tensor R \to L_Y$
is a surjection. 
By Theorem~\ref{thmPreciseA}, $L_Y \subset H_*(\lY;R)$ is torsion-free.
As a result, there is an induced surjection $\tilde{h}_Y: F(\pi_*(\lY)
\tensor R) \onto L_Y$.

Furthermore if we tensor with $\Q$ we see that $\pi_*(\lY) \tensor \Q
\isomto L_Y \tensor \Q$, which is a result of Cartan and Serre
(see~\cite{fhtRHT}).
Thus $F(\pi_*(\lY) \tensor R)$ and $L_Y$ have the same Hilbert series,
and therefore $\tilde{h}_Y$ is an isomorphism of Lie algebras.
 
So using Theorem~\ref{thmHurewicz}(b), the canonical algebra map
\[ U(F(\pi_*(\lY) \tensor R)) \to UL_Y \isomto H_*(\lY;R)
\] 
is an isomorphism.
\end{proof}

Corollaries \ref{corFiniteCellComplex} and \ref{corFelixAvramov}
follow immediately from Theorem~\ref{thmHurewicz}.

\section{Examples} \label{sectionEg}

\emph{(Spherical) $n$-cones} are those spaces $X$
such that there exists a finite sequence
\[
* = X_0 \subset X_1 \subset \cdots \subset X_n \simeq X 
\]
where for $k \geq 0$, $X_{k+1}$ is the \emph{adjunction space}
\[
X_{k+1} = X_k \cup_{f_{k+1}} \bigl( \bigvee_{j \in J} e^{n_j+1} \bigr).
\]
In particular, any finite CW-complex is an $n$-cone for some $n$.

\begin{eg} \label{egSemiInertTop} 
Let $X = S^3_a \vee S^3_b$ and let $\iota_a, \ \iota_b$ denote the
inclusions of the spheres into $X$.
Let $Y = X \cup_{\alpha_1 \vee \alpha_2} (e^8 \vee e^8)$ where the
attaching maps are given by the iterated Whitehead products $\alpha_1
= [[\iota_a, \iota_b], \iota_a]$ and $\alpha_2 = [[\iota_a,
\iota_b],\iota_b]$.
\end{eg}

Let $W = S^7 \vee S^7$ and for $i=1,2$ let $\hat{\alpha}_i: S^6 \to
\lX$ denote the adjoint of $\alpha_i$.
Let $R=\Fp$ with $p>3$ or $R \subset \Q$ containing $\frac{1}{6}$.
Then $[L^W_X] = [R\{h_X(\hat{\alpha}_1),h_X(\hat{\alpha}_2)\}]$.

$Y$ has an Adams-Hilton model (see Section~\ref{sectionAH}) $U(L,d)$
where $L = \freeL \langle x,y,a,b \rangle, \ |x|=|y|=2, \ dx=dy=0, \
da=[[x,y],x]$ and $db=[[x,y],y]$.
Furthermore $h_X(\hat{\alpha}_1) = [da]$ and
$h_X(\hat{\alpha}_2)=[db]$.

It is well-known that over a field, a Lie subalgebra of a free Lie
algebra is also free.
Thus, since $\freeL \langle x,y \rangle$ is a free Lie algebra, 
$UL$ is a free dga extension (in the sense of
Definition~\ref{defnFreeDGAextension}). 

Let $u = [a,y] - [b,x]$.
Then $du = [[[x,y],x],y] - [[[x,y],y],x] = [[x,y],[x,y]] = 0$.
Since $u$ is not a boundary $0 \neq [u] \in (HL)_1$ and $0 \neq [u]
\in (HUL)_1$.
Thus $UL$ is not an inert dga extension.
By the definition of homology 
\[
(HL)_0 \isom \freeL \langle x,y
\rangle \left/ \, \bigl[ R\{[[x,y],x],[[x,y],y]\} \bigr]. \right.
\]

One can check that $(HL)_1$ is freely generated by the $(HL)_0$-action
on $[u$].
Thus $UL$ is a semi-inert extension.
Therefore by Theorem~\ref{thmDGAextension}(b),
\[ 
HUL \isom U( (HL)_0 \amalg \freeL \langle [u] \rangle )
\]
as algebras.
Thus as algebras 
\begin{eqnarray*}
H_*(\lY;R) & \isom & U \left( L_X/[L^W_X] \amalg \freeL \langle [u]
\rangle \right) \\
& \isom & \left( H_*(\lX;R)/ (L^W_X) \right) \amalg \freeT \langle [u]
\rangle.
\end{eqnarray*}
$\alpha_1 \vee \alpha_2$ is a non-inert, semi-inert attaching map.

If $R = \Z[\frac{1}{6}]$ then since $X$ is a wedge of spheres there
exists a map $\sigma_X$ right inverse to $h_X$.
By Lemma~\ref{lemmaIp}, $P_Y = \{2,3\}$.
By Theorem~\ref{thmHurewicz}, 
\[ H_*(\lY;R) \isom UL_Y, \quad L_Y \isom L^X_Y \amalg \freeL \langle
{w} \rangle 
\isom \freeL \langle x,y,w \rangle / J,
\]
where $w = h_Y(\hat{\omega})$ with $\hat{\omega}$ the adjoint of some
map $\omega: S^{10} \to Y$,
and $J$ is the Lie ideal generated by $\{[[x,y],x], [[x,y],y]\}$.
Furthermore localized at R, $\lY \in \PiS$ and there exists a map
$\sigma_Y$ right inverse to $h_Y$.
\eolBox

\begin{eg} \label{egFatWedge}
The $6$-skeleton of $S^3 \times S^3 \times S^3$.
\end{eg}

This space $Y$ is also known as the \emph{fat wedge} $FW(S^3, S^3, S^3)$. 
Let $X$ = $S^3_a \vee S^3_b \vee S^3_c$.
Let $\iota_a, \iota_b$ and $\iota_c$ be the inclusions of the
respective spheres.
Let $W = \bigvee_{j=1}^3 S^5$.
Then $Y = X \cup_f (\bigvee_{j=1}^3 e^6_j)$ where $f: W \to X$ is given by
$\bigvee_{j=1}^3 \alpha_j$ with $\alpha_1 = [\iota_b,\iota_c]$, $\alpha_2 =
[\iota_c,\iota_a]$ and $\alpha_3 = [\iota_a,\iota_b]$.

Let $R = \Z[\frac{1}{6}]$.
Then $Y$ has Adams-Hilton model $U(L,d)$ where $L = \freeL \langle
x,y,z,a,b,c \rangle$,  $|x| = |y| = |z| = 2,  \ dx=dy=dz=0, \
da=[y,z], \ db=[z,x]$ and $dc=[x,y]$. 
Let $w = [x,a] + [y,b] + [z,c]$.

By the same argument as in the previous example, as algebras
\[ H_*(\lY;R) \isom UL_Y, \quad L_Y \isom L^X_Y \amalg \freeL \langle
w \rangle
\isom \freeL \langle x,y,z,w \rangle / J,
\]
where $w = h_Y(\hat{\omega})$ with $\hat{\omega}$ the adjoint of some
map $\omega: S^8 \to Y$ and  
$J$ is the Lie ideal generated by $\{[x,y],[y,z],[z,x]\}$.
Furthermore localized at R, $\lY \in \PiS$ and there exists a map
$\sigma_Y$ right inverse to $h_Y$.
\eolBox

The following \emph{spherical three-cone} $Y$, illustrates our results.

\begin{eg} \label{egSemiInert3coneConstruct}
Let $R = \Z[\frac{1}{6}]$.
All spaces here are localized at $R$.
For $i=1,2$ let $Z_i = A_i \cup_{\alpha_1 \vee \alpha_2} (e^8 \vee
e^8)$ be two copies of the two-cone from Example~\ref{egSemiInertTop}.
Let $X = Z_1 \vee Z_2$, $W = S^{28} \vee S^{28}$ and let $f = \beta_1
\vee \beta_2$ where $\beta_1 = [[\omega_1, \omega_2], \omega_1]$ and 
$\beta_2 = [[\omega_1, \omega_2], \omega_2]$. 
Let $Y = X \cup_f (e^{29} \vee e^{29})$.
\end{eg}

Now, 
\begin{equation} \label{eqnLXinSIeg} 
L_X \isom L_{Z_1} \amalg L_{Z_2} \isom L^{A_1}_{Z_1} \amalg
L^{A_2}_{Z_2} \amalg \freeL \langle w_1, w_2 \rangle.
\end{equation}
If follows from this that $f$ is a free attaching map.
Thus $Y$ satisfies the hypotheses of Theorem~\ref{thmPreciseA}.
Using Theorem~\ref{thmPreciseA} and Anick's formula 
%
%
one can calculate that if $f$ is semi-inert then
$K'(z) = z^{37}$.

Recall that $\eL = (L_X \amalg \freeL \langle e,g \rangle, d')$ where
$d' e = [[w_1,w_2],w_1]$ and $d' g = [[w_1,w_2],w_2]$.
Also recall (from Theorem~\ref{thmPreciseA}) that $(H\eL)_0 \isom L^X_Y$.
Let $\bar{u} = [e,w_2] + [g,w_1]$ (with $|\bar{u}| = 37$).
Then one can check that $d\bar{u}=0$ and
$[\bar{u}]$ is a basis for $(H\eL)_1$ as a free $(H\eL)_0$-module;
so $f$ is indeed semi-inert.
Let $\sigma_X = \sigma_{Z_1} \oplus \sigma_{Z_2}$.
It is right inverse to $h_X$.
By Lemma~\ref{lemmaIp}, $P_Y = \{2,3\}$.
As a result by Theorem~\ref{thmHurewicz},
\begin{gather*} 
H_*(\lY;R) \isom UL_Y \text{ where } \\
L_Y \isom L^X_Y \amalg \freeL \langle u \rangle
\isom \freeL \langle x_1, y_1, x_2, y_2, w_1, w_2, u \rangle / J,
\end{gather*}
with $u = h_Y(\hat{\mu})$ for some map $\mu: S^{38} \to Y$ and $J$ the
Lie ideal generated by  
\begin{multline*} 
\{ [[x_1,y_1],x_1], \ [[x_1,y_1],y_1], \ [[x_2,y_2],x_2], \\
[[x_2,y_2],y_2], \ [[w_1,w_2],w_1], \ [[w_1,w_2],w_2] \}.
\end{multline*}
Furthermore $\lY \in \PiS$ and there exists a map $\sigma_Y$ right
inverse to $h_Y$.  

Note that
\[
L_Y \isom L_1^1 \amalg L_1^2 \amalg L_2 \amalg \freeL \langle u
\rangle
\]
where $L_1^1 \isom L_1^2 \isom 
R\{ x,y,[x,y] \}$
and
$L_2 \isom 
R \{ w_1, w_2, [w_1,w_1], [w_1,w_2], [w_2,w_2] \}$.
\eolBox

\begin{eg} \label{egSemiInertNCones}
An infinite family of finite CW-complexes constructed out of
semi-inert attaching maps
\end{eg}

The construction in the previous example can be extended inductively. 
By induction, we will construct spaces $X_n$ and maps $\omega_n:
S^{\lambda_n} \to X_n$ for $n \geq 1$ such that $X_n$ is an $n$-cone
constructed out of a sequence of semi-inert attaching maps.
Given $\omega_n$, let $w_n = h_{X_n}([\omega_n])$ and given $w_i^a$ and
$w_i^b$, 
let $L_i = \freeL \langle w_i^a, w_i^b \rangle/ J_i$ where $J_i$ is
the Lie ideal of brackets of bracket length $\geq 3$.

Let $R = \Z[\frac{1}{6}]$.
Begin with $X_1 = S^3$ and $\lambda_1 = 3$.  
Let $\omega_1: S^{\lambda_1} \to X_1$ be the identity map.

Given $X_n$, let $X_n^a$ and $X_n^b$ be two copies of $X_n$.
For $n\geq 1$, let 
\[
X_{n+1} = X_n^a \vee X_n^b \cup_{f_{n+1}} \left(
e^{\kappa_{n+1}} \vee e^{\kappa_{n+1}} \right),
\]
where $\kappa_{n+1} = 3 \lambda_n - 1$ and $f_{n+1} = [[\omega_n^a,
\omega_n^b], \omega_n^a] \vee [[\omega_n^a, \omega_n^b], \omega_n^b]$.

By the same argument as in the previous example, 
$f_{n+1}$ is a semi-inert cell attachment and
there exists a map 
\[
\omega_{n+1}: S^{\lambda_{n+1}} \to X_{n+1}
\] 
where $\lambda_{n+1} = 4 \lambda_n - 2$, such that 
\[
H_*(\lX_{n+1}; R) \isom UL_{X_{n+1}} \text{ where } L_{X_{n+1}} =
\Bigl( \coprod_{ \substack{1 \leq i \leq n \\ 1 \leq j \leq 2^{n-i}} }
L^j_i \Bigr) 
\amalg \freeL \langle w_{n+1} \rangle
\]
with $L^j_i$ a copy of $L_i$ and $w_{n+1} = h_{X_{n+1}} (\hat{\omega}_{n+1})$.
\eolBox

\begin{eg}
An uncountable family of CW-complexes constructed out of semi-inert
cell attachments
\end{eg}


At each stage of the inductive construction in
Example~\ref{egSemiInertNCones}, we could have used an attachment of
the type in Example~\ref{egFatWedge} instead of an attachment of the
type in Example~\ref{egSemiInertTop}. 
For $\alpha \in [0,1)$ use the binary expansion of $\alpha$ to choose
the sequence of attachments to obtain a space $X_{\alpha}$.

\begin{eg}
CW-complexes $X$ with only odd-dimensional cells
\end{eg}

Let $X^{(n)}$ denote the $n$-skeleton of $X$.
From the CW-structure of $X$, there is a sequence of cell attachments for
$k\geq 1$. 
\[ W_{2k} \xto{f_{2k+1}} X^{(2k-1)} \to X^{(2k+1)}
\]
where $W_{2k}$ is a finite wedge of $(2k)$-dimensional spheres,
$X^{(2k+1)}$ is the adjunction space of the cell attachment and
$X^{(1)} = *$.  
Assume $R \subset \Q$ containing $\frac{1}{6}$.
Let $P_{X^{(n)}}$ be the set of implicit primes of $X^{(n)}$ (see
Definition~\ref{defnIp}). 
We will show by induction that
\[
H_*(\lX^{(2k+1)}; \Z[{P_{X^{(2k+1)}}}^{-1}]) \isom UL_{X^{(2k+1)}}
\]
where $L_{X^{(2k+1)}} \isom \freeL V^{(2k+1)}$ with 
$V^{(2k+1)}$ concentrated in even dimensions.
In addition localized away from $P_{X^{(2k+1)}}$, $\lX^{(2k+1)} \in \PiS$
and there exists a map $\sigma_{X^{(2k+1)}}$ right inverse to
$h_{X^{(2k+1)}}$.

For $k=0$ these conditions are trivial. 
Assume they hold for $k-1$.
Let $\eL = (L_{X^{(2k-1)}} \amalg \freeL K, d')$ where $K$ is a free
$R$-module in dimension $2k$ corresponding to the spheres in $W_{2k}$.
For degree reasons $L^{W_{2k}}_{X^{(2k-1)}} = d'(K) = 0$.
So $f_{2k+1}$ is automatically free.
Furthermore $\eL$ has zero differential so $H\eL = \freeL V^{(2k-1)}
\amalg \freeL K$.
Thus $f_{2k+1}$ is semi-inert. 
By Theorem~\ref{thmHurewicz},
$H_*(\lX^{(2k+1)}; \Z[{P_{X^{(2k+1)}}}^{-1}]) \isom UL_{X^{(2k+1)}}$, 
where 
\[
L_{X^{(2k+1)}} \isom L^{X^{(2k-1)}}_{X^{(2k+1)}} \amalg \freeL K
\isom L_{X^{(2k-1)}} \amalg \freeL K \isom \freeL (V^{(2k-1)} \oplus
K).
\]
Also by Theorem~\ref{thmHurewicz}, localized away from
$P_{X^{(2k+1)}}$, $\lX^{(2k+1)} \in \PiS$, and there exists a map
$\sigma_{X^{(2k+1)}}$ right inverse to  $h_{X^{(2k+1)}}$.

Therefore by induction $H_*(\lX;\Z[{P_X}^{-1}]) \isom UL_X$, where $L_X
\isom \freeL (s^{-1}\tilde{H}_*(X))$ with $s$ the suspension map,
and localized away from $P_X$, $\lX \in \PiS$.  
\eolBox

\bibliographystyle{amsalpha}
\bibliography{my}

\providecommand{\bysame}{\leavevmode\hbox to3em{\hrulefill}\thinspace}
\providecommand{\MR}{\relax\ifhmode\unskip\space\fi MR }
\providecommand{\MRhref}[2]{%
  \href{http://www.ams.org/mathscinet-getitem?mr=#1}{#2}
}
\providecommand{\href}[2]{#2}
\begin{thebibliography}{FHT01}

\bibitem[AH56]{adamsHilton}
J.~F. Adams and P.~J. Hilton, \emph{On the chain algebra of a loop space},
  Comment. Math. Helv. \textbf{30} (1956), 305--330. \MR{17,1119b}

\bibitem[Ani82]{anickThesis}
David~J. Anick, \emph{A counterexample to a conjecture of {S}erre}, Ann. of
  Math. (2) \textbf{115} (1982), no.~1, 1--33. \MR{86i:55011a}

\bibitem[Ani89]{anickCat2}
\bysame, \emph{Homotopy exponents for spaces of category two}, Algebraic
  topology (Arcata, CA, 1986), Lecture Notes in Math., vol. 1370, Springer,
  Berlin, 1989, pp.~24--52. \MR{90c:55010}

\bibitem[Ani92]{anickSLSD}
\bysame, \emph{Single loop space decompositions}, Trans. Amer. Math. Soc.
  \textbf{334} (1992), no.~2, 929--940. \MR{93g:55011}

\bibitem[Avr82]{avramovFreeLieSubalgebras}
Luchezar~L. Avramov, \emph{Free {L}ie subalgebras of the cohomology of local
  rings}, Trans. Amer. Math. Soc. \textbf{270} (1982), no.~2, 589--608.
  \MR{83g:13010}

\bibitem[Bau81]{bauesCommutator}
Hans~Joachim Baues, \emph{Commutator calculus and groups of homotopy classes},
  London Mathematical Society Lecture Note Series, vol.~50, Cambridge
  University Press, Cambridge, 1981. \MR{83b:55012}

\bibitem[Bub03]{bubenikThesis}
Peter Bubenik, \emph{Cell attachments and the homology of loop spaces and
  differential graded algebras}, Ph.D. thesis, University of Toronto, 2003.

\bibitem[FHT84]{fhtCat2}
Yves F{\'e}lix, Stephen Halperin, and Jean-Claude Thomas, \emph{Sur l'homotopie
  des espaces de cat\'egorie {$2$}}, Math. Scand. \textbf{55} (1984), no.~2,
  216--228. \MR{86k:55006}

\bibitem[FHT01]{fhtRHT}
\bysame, \emph{Rational homotopy theory}, Graduate Texts in Mathematics, vol.
  205, Springer-Verlag, New York, 2001. \MR{2002d:55014}

\bibitem[FT89]{felixThomasAttach}
Yves F{\'e}lix and Jean-Claude Thomas, \emph{Effet d'un attachement cellulaire
  dans l'homologie de l'espace des lacets}, Ann. Inst. Fourier (Grenoble)
  \textbf{39} (1989), no.~1, 207--224. \MR{90j:55012}

\bibitem[HL87]{halperinLemaireInert}
Stephen Halperin and Jean-Michel Lemaire, \emph{Suites inertes dans les
  alg\`ebres de {L}ie gradu\'ees (``{A}utopsie d'un meurtre. {II}'')}, Math.
  Scand. \textbf{61} (1987), no.~1, 39--67. \MR{89e:55022}

\bibitem[HL96]{hessLemaireNice}
Kathryn Hess and Jean-Michel Lemaire, \emph{Nice and lazy cell attachments}, J.
  Pure Appl. Algebra \textbf{112} (1996), no.~1, 29--39. \MR{97e:55006}

\bibitem[Lem78]{lemaireAutopsie}
Jean-Michel Lemaire, \emph{``{A}utopsie d'un meurtre'' dans l'homologie d'une
  alg\`ebre de cha\^\i nes}, Ann. Sci. \'Ecole Norm. Sup. (4) \textbf{11}
  (1978), no.~1, 93--100. \MR{58 \#18423}

\bibitem[MM65]{milnorMooreHopfAlgebras}
John~W. Milnor and John~C. Moore, \emph{On the structure of {H}opf algebras},
  Ann. of Math. (2) \textbf{81} (1965), 211--264. \MR{30 \#4259}

\bibitem[Sco02]{scottBSS}
Jonathan~A. Scott, \emph{Algebraic structure in the loop space homology
  {B}ockstein spectral sequence}, Trans. Amer. Math. Soc. \textbf{354} (2002),
  no.~8, 3075--3084 (electronic). \MR{2003c:55008}

\bibitem[Sco03]{scottTFMM}
\bysame, \emph{A torsion-free {M}ilnor-{M}oore theorem}, J. London Math. Soc.
  (2) \textbf{67} (2003), no.~3, 805--816. \MR{1 967 707}

\bibitem[Whi39]{whiteheadSimplicialSpaces}
J.~H.~C. Whitehead, \emph{Simplicial spaces, nuclei and {$m$}-groups}, Proc.
  London Math. Soc. (2) \textbf{45} (1939), 243--327.

\bibitem[Whi41]{whiteheadAddingRelations}
\bysame, \emph{On adding relations to homotopy groups}, Ann. of Math. (2)
  \textbf{42} (1941), 409--428. \MR{2,323c}

\end{thebibliography}

\end{document}